\documentclass{amsart}
\usepackage{graphicx}
\newtheorem{thm}{Theorem}
\newtheorem{lem}{Lemma}

\newtheorem{prop}{Proposition}
\newtheorem{cor}{Corollary}

\newtheorem{conj}{Conjecture}
\renewcommand{\thethm}{\Alph{thm}}
\begin{document}
\setlength{\textheight}{7.7truein}  

\title{Superbridge index of composite knots}
\author[Gyo Taek Jin]{Gyo Taek Jin\\
\vspace{0.7truein}}
\address{
\parbox{0.8\linewidth}{Department of Mathematics\\
Korea Advanced Institute of Science and Technology\\
Taejon, 305-701, Korea\\}}
\email{trefoil@kaist.ac.kr}

\begin{abstract}
An upper bound of the superbridge index of the connected sum of two knots
is given in terms of the braid index of the summands. Using this upper bound
and minimal polygonal presentations, we give an upper bound in terms of the
superbridge index and the bridge index of the summands when they are torus
knots.
In contrast to the fact that the difference between the sum of bridge indices
of two knots and the bridge index of their connected sum is always one,
the corresponding difference for the superbridge index can be arbitrarily large.
\end{abstract}
\maketitle
\thispagestyle{empty}
\section{Introduction}
\label{intro}
Throughout this article a {\em knot\/} is a piecewise smooth simple closed
curve embedded in the three dimensional Euclidean space $\mathbb R^3$.
Two knots are {\em equivalent\/} if there is a
piecewise smooth autohomeomorphism of $\mathbb R^3$ mapping one knot
onto the other. The equivalence class of a knot $K$ will be called the
{\em knot type\/} of $K$ and denoted by~$[K]$.

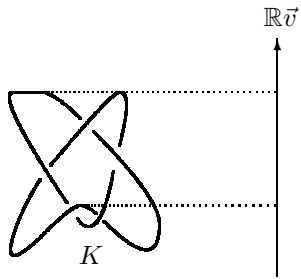
\begin{figure}[t]
\begin{picture}(105,100)(20,40)
\put(120,40){\vector(0,1){90}}
\put(115,135){$\mathbb R\vec v$}
{\thicklines
\qbezier(40,68)(12,110)(20,110)
\qbezier(20,110)(25,110)(30,110)
\qbezier(30,110)(35,110)(45,100)
\qbezier(50,95)(75,70)(75,60)
\qbezier(75,60)(75,40)(54,61)
\qbezier(51,64)(48,67)(45,67)
\qbezier(45,67)(42,67)(32,57)
\qbezier(32,57)(12,37)(22,62)
\qbezier(22,62)(26,71)(30,77)
\qbezier(34,82.4)(57,110)(60,110)
\qbezier(60,110)(65,110)(60,90)
\qbezier(57.5,80)(54,59)(48,59)
\qbezier(48,59)(46,59)(44,62)
}
\qbezier[45](30,110)(65,110)(119,110)
\qbezier[38](45,67)(65,67)(119,67)
\put(45,45){$K$}
\end{picture}
\caption{$b_{\vec v}(K)=3$}\label{fig:local-max}
\end{figure}

The {\em crookedness\/} of a knot $K$ embedded in $\mathbb R^3$ with
respect to a unit vector $\vec v$ is the number of connected components of
the preimage of the set of local maximum values of the orthogonal projection
$K\to \mathbb R\vec v$, denoted by $b_{\vec v}(K)$.
Figure~\ref{fig:local-max} illustrates an example.
For any open subarc $S$ of a knot $K$, the crookedness of $S$
with respect to $\vec v$, denoted by $b_{\vec v}(K\mid S)$,
can be defined similarly using the projection $S\to \mathbb R\vec v$.
The {\em superbridge number\/} and the {\em superbridge index\/} of $K$,
denoted by $s(K)$ and $s[K]$, are defined to be
``$\max b_{\vec v}(K)$'' and ``$\min\,\max b_{\vec v}(K)$'', respectively,
where the maximum is taken over all unit vectors and the minimum taken over
all equivalent embeddings of $K$.
This invariant was introduced by Kuiper~\cite{kuiper} who computed the
superbridge index for all torus knots.

\begin{thm}[Kuiper]\label{thm:kuiper-torus}
For any two coprime integers $p$ and $q$, satisfying $2\le p<q$,
the torus knot of type $(p,q)$ has superbridge index $\min\{2p, q\}$.
\end{thm}

The {\em bridge index\/} $b[K]$ can be defined in a similar way by
``$\min\,\min b_{\vec v}(K)$.''
One of the most well-known theorem about bridge index is

\begin{thm}[Schubert]\label{thm:schubert}
Given two knots $K_1$ and $K_2$, any connected
sum\footnote{Since $K_1$ and $K_2$ are not oriented, their (unoriented)
connected sum may not be unique.} $K_1\sharp K_2$ satisfies
$$b[K_1\sharp K_2]=b[K_1]+b[K_2]-1.$$
\end{thm}

This work is an attempt to find a similar formula for the
superbridge index.
A proof of Schubert's theorem in a more generalized context can be found
in~\cite{doll}.

\medskip
Let $\beta[K]$ denote the {\em braid index}, i.e., the minimal number of
strings among all braids whose closures are equivalent to $K$. According to
Kuiper, the superbridge index of a nontrivial knot is always greater than the
bridge index and not greater than twice the braid index~\cite{kuiper}.

\begin{thm}[Kuiper]\label{thm:kuiper-braid} If $K$ is a nontrivial knot, then
$$b[K]<s[K]\le2\beta[K].$$
\end{thm}

Kuiper used Milnor's total curvature to prove the first
inequality~\cite{milnor1}.
The closed braid constructed by Kuiper used to prove the second inequality is
discussed in Section~\ref{sec:closed-braids}.
From Theorem~\ref{thm:schubert} and Theorem~\ref{thm:kuiper-braid}, we obtain

\begin{cor}\label{cor:at-least-4}
If $K_1$ and $K_2$ are nontrivial knots, any connected sum $K_1\sharp K_2$
satisfies the inequality
$$s[K_1\sharp K_2]\ge4.$$
\end{cor}

\section{Theorems and Conjectures}
\renewcommand{\thethm}{\arabic{thm}}
\setcounter{thm}{0}
\begin{thm}\label{thm:braid-sum}
If $K_1$ and $K_2$ are nontrivial knots,
any connected sum $K_1\sharp K_2$ satisfies the
inequality
$$s[K_1\sharp K_2]\le
\max\{2\beta[K_1]+\beta[K_2],\beta[K_1]+2\beta[K_2]\}-1.$$
\end{thm}

\begin{thm}\label{thm:torus-sum}
If $K_1, K_2$ are torus knots, any connected sum $K_1\sharp K_2$ satisfies
the inequality
$$s[K_1\sharp K_2]\le\max\{s[K_1]+b[K_2],b[K_1]+s[K_2]\}-1.$$
\end{thm}
The next corollary shows that the equality in Theorem~\ref{thm:torus-sum}
holds in infinitely many cases.
\begin{cor}\label{cor:sb=b+1}
Let $p_i\ge2$ and let $K_i$ be the torus knot of type $(p_i,p_i+1)$,
for $i=1,2$. Then
$$s[K_1\sharp K_2]=p_1+p_2.$$
\end{cor}

\noindent{\bf Proof:} By Theorem~\ref{thm:kuiper-torus}, $s[K_i]=p_i+1$.
Since $b[K_i]=p_i$, from Theorem~\ref{thm:schubert}, Theorem~\ref{thm:kuiper-braid} and Theorem~\ref{thm:torus-sum}, we obtain
$p_1+p_2-1<s[K_1\sharp K_2]\le p_1+p_2$.\qed

\medskip
Using the first inequality in Theorem~\ref{thm:kuiper-braid},
we obtain the following generalization of \cite[Corollary~11]{jin-poly}.

\begin{cor}\label{cor:torus-sum}
If $K_1, K_2$ are torus knots, any connected sum $K_1\sharp K_2$ satisfies
the inequality
$$s[K_1\sharp K_2]\le s[K_1]+s[K_2]-2.$$
\end{cor}

\medskip
The inequality in Theorem~\ref{thm:torus-sum} is equivalent to
$$s[K_1]+s[K_2]-s[K_1\sharp K_2]\ge\min\{s[K_1]-b[K_1],s[K_2]-b[K_2]\}+1.$$
If $K_i$ is a torus knot of type $(p_i,q_i)$ with $2\le p_i<q_i$,
the right hand side of the above inequality is equal to
$\min\{p_1,p_2,q_1-p_1,q_2-p_2\}+1$,
which can be arbitrarily large. Therefore we have
\begin{cor}
The difference $s[K_1]+s[K_2]-s[K_1\sharp K_2]$
can be arbitrarily large.
\end{cor}
We conjecture that Theorem~\ref{thm:torus-sum} and
Corollary~\ref{cor:torus-sum} are true for any knots:
\begin{conj}\label{conj:connected-sum-harder} Any connected sum
of two knots $K_1$ and $K_2$ satisfies the inequality
$$s[K_1\sharp K_2]\le\max\{ s[K_1]+b[K_2],b[K_1]+s[K_2]\}-1.$$
\end{conj}
\begin{conj}\label{conj:connected-sum-easier} If $K_1$ and $K_2$ are
nontrivial knots, any connected sum $K_1\sharp K_2$
satisfies the inequality
$$s[K_1\sharp K_2]\le s[K_1]+s[K_2]-2.$$
\end{conj}

As Corollary~\ref{cor:torus-sum} follows from Theorem~\ref{thm:torus-sum},
Conjecture~\ref{conj:connected-sum-easier} follows from
Conjecture~\ref{conj:connected-sum-harder}.
The readers may wonder whether the inequality
$$s[K_1\sharp K_2]\ge\max\{ s[K_1]+b[K_2],b[K_1]+s[K_2]\}-1$$
would be true. So far no reasonable lower bound formula for $s[K_1\sharp K_2]$
has been found. We do not even know if the following is true.

\begin{conj}\label{conj:connected-sum-lower} If $K_1$ and $K_2$ are
nontrivial knots, any connected sum $K_1\sharp K_2$
satisfies the inequality
$$s[K_1\sharp K_2]> \max\{s[K_1],s[K_2]\}.$$
\end{conj}

\begin{table}[t]
\def\d{\relax}
\def\e{\relax}
\def\s{$\star$}
\def\c{Corollary~\ref{cor:torus-sum}}
\def\t{Theorem~\ref{thm:torus-sum}}
\begin{tabular}{|l|c|c|c|c|}
\hline
factors of $K$&$s[K]$&lower bound&upper bound \\
\hline
$3_1\quad 3_1$&4\d&$b[K]=3$&\c\\
$3_1\quad 4_1$&4\e&$b[K]=3$&$p[K]=9$\\
$3_1\quad 5_1$& \phantom\s5\s&$s[5_1]=4$&\t\\
$3_1\quad 7_1$& \phantom\s5\s&$s[7_1]=4$&\t\\
$3_1\quad 7_6$& \phantom\s5\s&$s[7_6]=4$&$p[K]\le11$\\
$3_1\quad 7_7$& \phantom\s5\s&$s[7_7]=4$&$p[K]\le11$\\
$3_1\quad 8_{16}$&5\e&$b[K]=4$&$p[K]\le11$\\
$3_1\quad 8_{17}$&5\e&$b[K]=4$&$p[K]\le11$\\
$3_1\quad 8_{18}$&5\e&$b[K]=4$&$p[K]\le11$\\
$3_1\quad 8_{19}$&5\d&$b[K]=4$&\c\\
$3_1\quad 8_{20}$&5\e&$b[K]=4$&$p[K]\le10$\\
$3_1\quad 8_{21}$&5\e&$b[K]=4$&$p[K]\le11$\\
$3_1\quad 9_{1}$&\phantom\s 5\s&$s[9_1]=4$&\t\\
$3_1\quad 9_{40}$&5\e&$b[K]=4$&$p[K]\le11$\\
$3_1\quad 9_{41}$&5\e&$b[K]=4$&$p[K]\le11$\\
$3_1\quad 9_{44}$&5\e&$b[K]=4$&$p[K]\le11$\\
$3_1\quad 9_{46}$&5\e&$b[K]=4$&$p[K]\le11$\\
$4_1\quad 5_1$&\phantom\s 5\s&$s[5_1]=4$&$p[K]\le11$\\
$4_1\quad 8_{19}$&5\e&$b[K]=4$&$p[K]\le11$\\
$4_1\quad 8_{20}$&5\e&$b[K]=4$&$p[K]\le11$\\
$5_1\quad 5_1$&\phantom\s5\s&$s[5_1]=4$&\t\\
$5_1\quad 7_1$&\phantom\s5\s&$s[7_1]=4$&\t\\
$7_1\quad 7_1$&\phantom\s5\s&$s[7_1]=4$&\t\\
$8_{19}\quad 8_{19}$&6\d&$b[K]=5$&\c\\
$3_1\quad 3_1\quad 3_1$&5\e&$b[K]=4$&$p[K]\le10$\\
$3_1\quad 3_1\quad 4_1$&5\e&$b[K]=4$&$p[K]\le11$\\
\hline
\end{tabular}
\medskip
\caption{}\label{tab:supbr}
\end{table}
\begin{figure}[t]
\includegraphics[bb=0 0 352 409,scale=0.35]{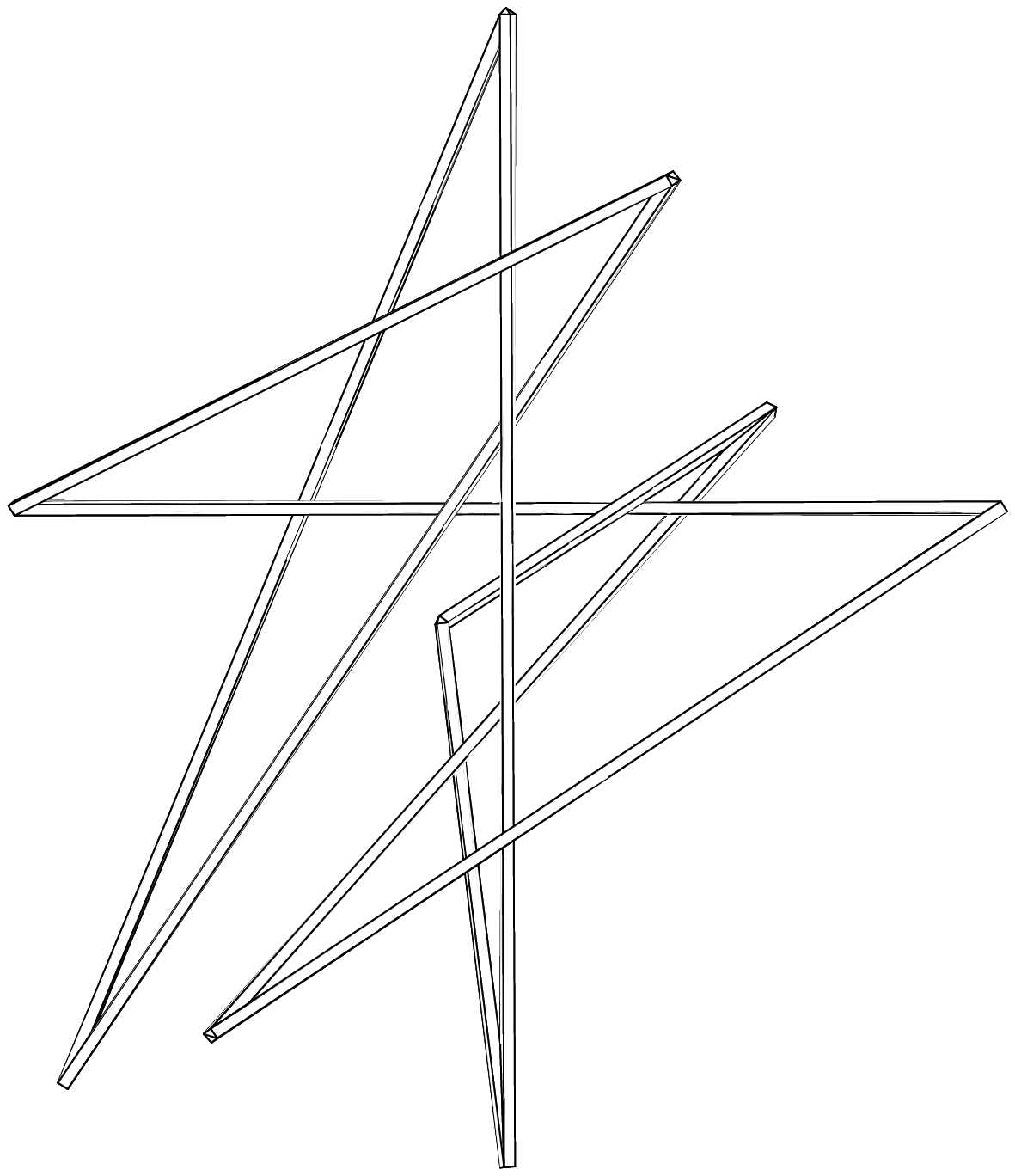}
\caption{A minimal polygonal connected sum of $3_1$ and $4_1$.}%
\label{fig:tref+fig8}
\end{figure}

\medskip
In Table~\ref{tab:supbr}, the symbols used for factors of $K$
indicate the prime knots as in
the knot tables of \cite{burde-zieschang,rolfsen}.
The knots $3_1,5_1,7_1,8_{19},9_1$ are torus knots of type
$(2,3), (2,5), (2,7), (3,4), (2,9)$, respectively.
Theorem~\ref{thm:torus-sum} is used to find upper bounds of superbridge index
for the connected sums of pairs of these knots.
There are three among them for which Corollary~\ref{cor:torus-sum}
also applies.  For the others, we used the inequality
\begin{equation}\label{eq:2s<p} 2\,s[K]\le p[K] \end{equation}
to find upper bounds,
where $p[K]$ is the {\em polygon index~}\cite{jin-poly,jin-kim}, i.e.,
the minimal
number of straight edges required to present the knot type of $K$.
Using the polygonal knots given in \cite{meissen,millett,scharein},
we verified that the inequality
$$p[K_1\sharp K_2]\le p[K_1]+p[K_2]-4$$
of \cite[Theorem~8]{jin-poly} can be applied to find upper bounds of $p[K]$
as given in the table.
The nine-edged polygonal knot\footnote{
It has vertices at
$(-30,0,-10)$, $(10,20,30)$, $(-27,-35,-70)$, $(0,30,10)$,
$(0,-40,10)$, $(-4,-7,8)$, $(16,6,-21)$, $(-18,-32,36)$,
$(30,0,-10)$. Figure~\ref{fig:tref+fig8} is its projection into the
$xy$-plane.}
of Figure~\ref{fig:tref+fig8}
is a connected sum of a trefoil knot and a figure eight knot.
It has polygon index 9 because it does not appear in the list of
\cite{calvo} containing all eight-edged knots.

The values marked with $\star$ are conjectured using
Theorem~\ref{thm:kuiper-torus}, Conjecture~\ref{conj:connected-sum-lower}
and \cite[Table~1]{jeon-jin-3}.
If Conjecture~\ref{conj:connected-sum-lower} is not true for any of them,
the correct value will be one less than as given in the table.
For all others, Theorem~\ref{thm:schubert} and Theorem~\ref{thm:kuiper-braid}
are used to determine strict lower bounds.

\medskip
The next two sections describe the constructions and their properties
required to prove
Theorem~\ref{thm:braid-sum} and Theorem~\ref{thm:torus-sum}.
Section~\ref{sec:proofs} contains the proofs.

\section{Closed braids}\label{sec:closed-braids}
\begin{figure}[t]
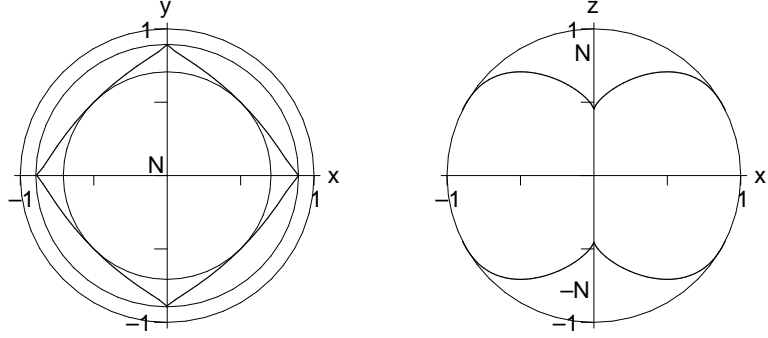

\includegraphics[242pt,335pt][364pt,456pt]{N1.ps}
\hskip0.5in
\includegraphics[242pt,335pt][364pt,456pt]{N2.ps}
\caption{Projections of $N$ and $-N$}\label{fig:boundary-N}
\end{figure}

Let {\bf i}, {\bf j}, {\bf k} denote the standard basis vectors of
$\mathbb R^3$ and let $\eta$ be the trivial knot given by the embedding
$(x,y)\mapsto(x,y,x^2)$ of the circle $x^2+y^2=1$. By \cite[Lemma~4.1]{kuiper},
we know that $s(\eta)=2$. Therefore, for any unit vector $\vec v$,
either $b_{\vec v}(\eta)=1$ or
$b_{\vec v}(\eta)=2$.
Let
$N=\{\vec v=v_1{\bf i}+v_2{\bf j}+v_3{\bf k}\in S^2\mid
      v_3>0,b_{\vec v}(\eta)=2\}.$
This is an open subset of $S^2$ satisfying the condition that
$b_{\vec v}(\eta)=2$ if and only if $\vec v\in N\cup(-N)$.
Two projections of $N$ and $-N$ are shown in Figure~\ref{fig:boundary-N}.

\begin{lem}\label{lem:boundary-N}
Let $G_{\rho,\alpha}(t)=-\rho\sin(t-\alpha)-(1-\rho^2)^{1/2}\sin2t$.
\begin{enumerate}
\item For any $\alpha$, there is a unique positive number
$\xi(\alpha)\in\left[1/\sqrt2,2/\sqrt5\,\right]$ such that the function
$G_{\xi(\alpha),\alpha}(t)$
has a multiple root.
\item $\partial N$ has a parametrization
$\alpha\mapsto\left(\xi(\alpha)\cos\alpha,
                     \xi(\alpha)\sin\alpha,
                     (1-\xi(\alpha)^2)^{1/2}\right)$.
\end{enumerate}
\end{lem}

\noindent{\bf Proof:} (a)
If $t_0$ is a multiple root of $G_{\rho,\alpha}(t)$, then
\begin{eqnarray*}
G_{\rho,\alpha}(t_0)&=&
-\rho\sin(t_0-\alpha)-(1-\rho^2)^{1/2}\sin2t_0=0,\\
G^\prime_{\rho,\alpha}(t_0)&=&
-\rho\cos(t_0-\alpha)-2(1-\rho^2)^{1/2}\cos2t_0=0.
\end{eqnarray*}
Eliminating $\alpha$, we get
$\rho=((1+3\cos^22t_0)/(2+3\cos^22t_0))^{1/2}$.
Therefore the inequality $1/\sqrt2\le\rho\le2/\sqrt5$ holds.

\begin{figure}[t]
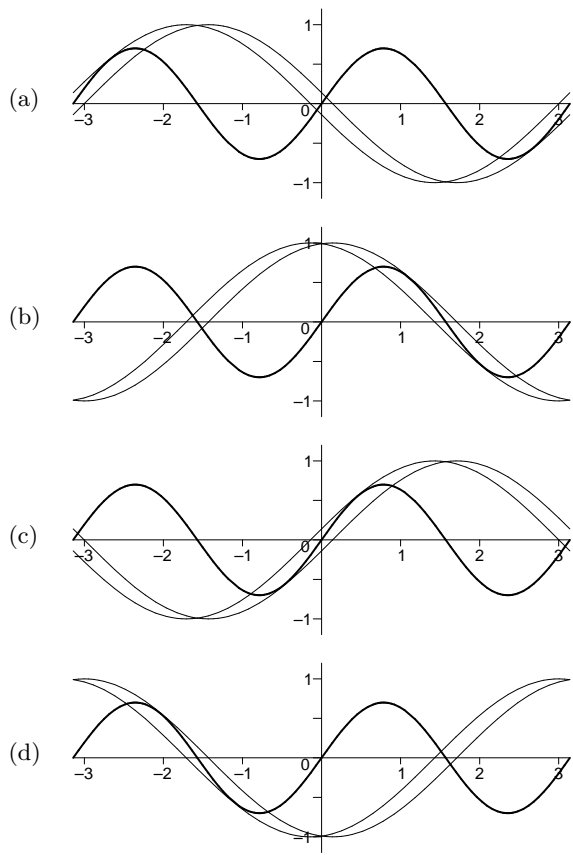

\small
\def\PutSines#1{%
\includegraphics[bb=168 340 446 465,scale=0.65]{#1}
}
\raisebox{1.25cm}{\hbox to 1cm{(a)\hfill}}\PutSines{2max-1.ps}

\raisebox{1.25cm}{\hbox to 1cm{(b)\hfill}}\PutSines{2max-2.ps}

\raisebox{1.25cm}{\hbox to 1cm{(c)\hfill}}\PutSines{2max-3.ps}

\raisebox{1.25cm}{\hbox to 1cm{(d)\hfill}}\PutSines{2max-4.ps}
\caption{$p(t)$'s and $q(t)$ with $1/2<(1/\rho^2-1)^{1/2}<1$}%
\label{fig:tangent-pairs}
\end{figure}

Suppose $1/\sqrt2<\rho<2/\sqrt5$, then $1/2<(1/\rho^2-1)^{1/2}<1$.
As illustrated in Figure~\ref{fig:tangent-pairs}, there are eight distinct
values of $\alpha$ modulo $2\pi$, such that the graphs of
$p(t)=-\sin(t-\alpha)$ and $q(t)=(1/\rho^2-1)^{1/2}\sin2t$
are tangent at some point. For these values of $\alpha$, the function
$G_{\rho,\alpha}(t)$ has double roots.

If $\alpha=k\pi\pm\pi/4$, $k=0,1$, then $\rho=1/\sqrt2$. In these cases,
the graphs of $p(t)$ and $q(t)$ are tangent at $t_0=\pi-\alpha$,
where $G_{\rho,\alpha}(t)$ has a double root.
If $\alpha=k\pi/2$, $k=0,1,2,3$, then $\rho=2/\sqrt5$. In these cases,
the graphs of $p(t)$ and $q(t)$ are tangent at $t_0=\pi-\alpha$,
where $G_{\rho,\alpha}(t)$ has a triple root.
This finishes the proof of part~(a) except the uniqueness which we omit.

\smallskip
(b) For a unit vector $\vec v=v_1{\bf i}+v_2{\bf j}+v_3{\bf k}$,
the projection $\eta\to\mathbb R\vec v$ is parametrized by
\begin{equation}\label{eqn:eta-by-vi}
f_{\vec v}(t)=v_1\cos t+v_2\sin t+v_3\cos^2t.
\end{equation}
Suppose $0<v_3<1$, then there is a unique number $\alpha_{\vec v}$
modulo $2\pi$ such that
$\cos\alpha_{\vec v}=v_1(1-v_3^2)^{-1/2}$ and
$\sin\alpha_{\vec v}=v_2(1-v_3^2)^{-1/2}$.
Substituting $\rho=(1-v_3^2)^{1/2}$, we get
\begin{equation}\label{eqn:eta-by-rho}
f_{\vec v}(t)=\rho\cos(t-\alpha_{\vec v})+(1-\rho^2)^{1/2}\cos^2t.
\end{equation}
If $\vec v\in\partial N$, then $f_{\vec v}^\prime(t)=0$ has a multiple root.
Since $f_{\vec v}^\prime(t)=G_{\rho,\alpha_{\vec v}}(t)$,
we know that $\partial N$ has the required parametrization.
The projection of $\partial N$ into the $xy$-plane in
Figure~\ref{fig:boundary-N} is the graph of the polar equation
$\rho=\xi(\alpha)$.\qed


\begin{lem}\label{lem:1-or-0}
Let $\vec v=v_1{\bf i}+v_2{\bf j}+v_3{\bf k}$ be a unit vector.
Then
$$b_{\vec v}(\eta\mid\eta_+)=
\begin{cases}
1&\text{if\/\ }\vec v\in N\text{\ or\/\ }\min\{v_1,v_3\}>0\\
0&\text{if\/\ }\vec v\notin N,v_1<0\text{\ and\/\ }v_3>0,
\end{cases}$$
where $\eta_+=\eta\cap\{(x,y,z)\mid x>0\}$.
\end{lem}

\medskip
\noindent{\bf Proof:} Again we use the parametrizations (\ref{eqn:eta-by-vi})
and (\ref{eqn:eta-by-rho}) for $\eta$. We have
\begin{eqnarray*}
f^\prime_{\vec v}(t)
&=&-v_1\sin t+v_2\cos t-v_3\sin2t\\
&=&-\rho\sin(t-\alpha_{\vec v})-(1-\rho^2)^{1/2}\sin2t.
\end{eqnarray*}

\noindent{\sc Case 1.} Suppose $\vec v\notin N$ and $v_1>0$.  Then
$f^\prime_{\vec v}(\pi/2)=-v_1<0<v_1=f^\prime_{\vec v}(-\pi/2)$.
Therefore $b_{\vec v}(\eta\mid\eta_+)=1$.

\noindent{\sc Case 2.} Suppose $v_3>1/\sqrt2$.
Since $0\le\rho<(1-\rho^2)^{1/2}$, we have
\begin{equation}\label{eqn:kpi-pm-pi/4}
f^\prime_{\vec v}(\pi/4)<0<f^\prime_{\vec v}(-\pi/4)
\text{\ and\ }
f^\prime_{\vec v}(5\pi/4)<0<f^\prime_{\vec v}(3\pi/4).
\end{equation}
Therefore there are two local maximum points, one in each of the two intervals
$(-\pi/4,\pi/4)$ and  $(3\pi/4,5\pi/4)$.
Therefore $\vec v\in N$ and $b_{\vec v}(\eta\mid\eta_+)=1$.

\noindent{\sc Case 3.} Suppose that $\vec v\in N$ and $v_3=1/\sqrt2$, then
$\rho=(1-\rho^2)^{1/2}=1/\sqrt2$ and $\alpha_{\vec v}\ne k\pi/2+\pi/4$
for any integer $k$. Therefore condition~(\ref{eqn:kpi-pm-pi/4}) holds,
and again we have $b_{\vec v}(\eta\mid\eta_+)=1$.

\begin{figure}[b]
\includegraphics[242pt,335pt][364pt,456pt]{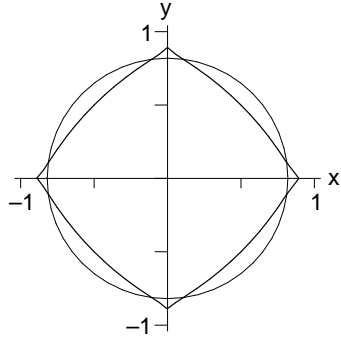}
\caption{$\partial N$ and the circle $x^2+y^2=\rho^2$}%
        \label{fig:eight-alphas}
\end{figure}

\noindent{\sc Case 4.}
Suppose that $\vec v\in N$ and $v_3<1/\sqrt2$. Then
$1/\sqrt5<v_3<1/\sqrt2$, hence $1/\sqrt2<\rho<2/\sqrt5$ and
$1/2<(1/\rho^2-1)^{1/2}<1$. The circle $x^2+y^2=\rho^2$ on the unit sphere
meets $\partial N$ at eight distinct points as
shown in Figure~\ref{fig:eight-alphas}.
Let $\alpha_0$ be the smallest positive number that $G_{\rho,\alpha_0}(t)$ has
double roots. Since $\vec v\in N$, it is on one of the four open
arcs of the circle inside $N$. These arcs correspond to the four
intervals for $\alpha_{\vec v}$ given in the table below.
\begin{center}
\def\s{\small}
\begin{tabular}{|c|c|c|c|}
\hline
\s\hbox to 0.9in{\hfill(a)\hfill}&\s\hbox to 0.9in{\hfill(b)\hfill}&
\s\hbox to 0.9in{\hfill(c)\hfill}&\s\hbox to 0.9in{\hfill(d)\hfill}\\
\hline
\s$|\alpha_{\vec v}|<\alpha_0$&
\s$|\alpha_{\vec v}-\pi/2|<\alpha_0$&
\s$|\alpha_{\vec v}-\pi|<\alpha_0$&
\s$|\alpha_{\vec v}-3\pi/2|<\alpha_0$\\
\hline
\end{tabular}
\end{center}
The four pairs of $p(t)$'s in Figure~\ref{fig:tangent-pairs} correspond to
the endpoints of these intervals. From Figure~\ref{fig:tangent-pairs},
we easily see that the sign of $f^\prime_{\vec v}(t)=\rho(p(t)-q(t))$
changes from positive to negative once in each of the intervals
$(-\pi/2,\pi/2)$ and $(\pi/2,3\pi/2)$. Therefore
$b_{\vec v}(\eta\mid\eta_+)=1$.

\noindent{\sc Case 5.}
Suppose $\vec v\notin N$ and $v_1\le0$. If $v_1=0$, any local extremum of
$f_{\vec v}$ occurs only at $(0,1,0)$ or $(0,-1,0)$.
If $v_1<0$, then
$f^\prime_{\vec v}(3\pi/2)=v_1<0<-v_1=f^\prime_{\vec v}(\pi/2)$.
Therefore $b_{\vec v}(\eta\mid\eta_-)=1$ where
$\eta_-=\eta\cap\{(x,y,z)\mid x<0\}$. Since $b_{\vec v}(\eta)=1$,
we obtain $b_{\vec v}(\eta\mid\eta_+)=0$.\qed

\medskip
Suppose $n$ is a positive integer and $K$ is a knot parametrized by
$$K(t)=((1+\lambda_1(t))\cos nt,(1+\lambda_1(t))\sin nt,\lambda_2(t)+\cos^2nt)$$
over any interval of length $2\pi$, for some smooth periodic functions
$\lambda_1$ and $\lambda_2$ with period $2\pi$ satisfying the conditions
\begin{eqnarray}
&&\lambda_1(t)^2+\lambda_2(t)^2<1\label{eqn:in-torus},\\
&&\lambda_1(t)=\lambda_2(t)=0\
               \text{if\ } |t|\le {3\pi}/{4n}\label{eqn:as-eta}, \\
&&\lambda_1(t),\ \lambda_2(t)\
        \text{are locally constant and negative\ } \nonumber\\
&&\phantom{\lambda_1(t),\ \lambda_2(t)\ } \label{eqn:like-eta}
       \text{if\ }{5\pi}/{4n}\le|t|\le\pi\ \text{and\ } \cos nt\ge-1/{\sqrt2}.
\end{eqnarray}
For any $\varepsilon$ with $0\le\varepsilon\le1$, we define
\begin{equation}
K^\varepsilon(t)=((1+\varepsilon\lambda_1(t))\cos nt,(1+\varepsilon\lambda_1(t))\sin nt,\varepsilon\lambda_2(t)+\cos^2nt).\label{eqn:Ke}
\end{equation}
Then $K^\varepsilon$ is a knot isotopic to $K$ and is the closure of
the $n$-braid $K^\varepsilon\cap\{(x,y,z)\mid x\le y\le -x\}$
when $0<\varepsilon\le1$. When $\varepsilon=0$, $K^\varepsilon$ is
an $n$-fold covering of $\eta$.
Since $K^\varepsilon_+=K^\varepsilon\cap\{(x,y,z)\mid x>0\}$
is the union of
$n$ disjoint parallel copies of $\eta_+$
up to radial scaling about the $z$-axis, we have
$b_{\vec v}(K^\varepsilon\mid K^\varepsilon_+)=
n\cdot b_{\vec v}(\eta\mid\eta_+),$
hence by Lemma~\ref{lem:1-or-0}, we obtain
\begin{equation}\label{eqn:n-or-0}
b_{\vec v}(K^\varepsilon\mid K^\varepsilon_+)=\begin{cases}
n&\text{if\ }\vec v\in N \text{\ or\ } \min\{v_1,v_3\}>0\\
0&\text{if\ }\vec v\notin N, v_1<0 \text{\ and\ }v_3>0
\end{cases}
\end{equation}
for any unit vector $\vec v=v_1{\bf i}+v_2{\bf j}+v_3{\bf k}$.

By \cite{kuiper}, we know that there is a
number $\varepsilon^\prime>0$ such that
$s(K^\varepsilon)=2n$ whenever
$0<\varepsilon\le\varepsilon^\prime$.
Let $0<\varepsilon\le\varepsilon^\prime$ and let
\begin{equation}\nonumber
N^\varepsilon=\{\vec v=v_1{\bf i}+v_2{\bf j}+v_3{\bf k}\in S^2\mid
      v_3>0,b_{\vec v}(K^\varepsilon)=2n\}.
\end{equation}
For any $\varepsilon$, $N^\varepsilon$ is an open set intersecting $N$ in
a neighborhood of $\bf k$. Since $N^0=N$ and is connected, $N^\varepsilon$
is also connected whenever $0<\varepsilon\le \varepsilon^{\prime\prime}$
for some $\varepsilon^{\prime\prime}\in(0,\varepsilon^\prime]$.

Suppose
$N^\varepsilon\cap\{(x,y,z)\mid x<0\}
\not\subset N\cap\{(x,y,z)\mid x<0\}$. Then there exists
a unit vector
$\vec v\in\partial N\cap N^\varepsilon\cap\{(x,y,z)\mid x<0\}$.
Since the projection $K^\varepsilon\to\mathbb R\vec v$ assumes no local
maximum in  $K^\varepsilon\cap\{(x,y,z)\mid x=0\}$,
and by the equation~(\ref{eqn:n-or-0}), we have
$$b_{\vec v}(K^\varepsilon\mid K^\varepsilon_-)
=b_{\vec v}(K^\varepsilon\mid K^\varepsilon_+)
+b_{\vec v}(K^\varepsilon\mid K^\varepsilon_-)
=b_{\vec v}(K^\varepsilon)=2n.
$$
There exists an open neighborhood $V$ of $\vec v$ contained in
$N^\varepsilon\cap\{(x,y,z)\mid x<0\}$ such that
\begin{equation}\label{eqn:2n-in-x<0}
b_{\vec u}(K^\varepsilon\mid K^\varepsilon_-)=2n
\end{equation}
for any $\vec u\in V$. For any $\vec u\in V\cap N$, we obtain the following
contradiction from (\ref{eqn:n-or-0}) and (\ref{eqn:2n-in-x<0}):
$$2n=b_{\vec u}(K^\varepsilon)
=b_{\vec u}(K^\varepsilon\mid K^\varepsilon_+)
+b_{\vec u}(K^\varepsilon\mid K^\varepsilon_-)
=3n.$$

\medskip
\begin{prop}\label{prop:Ne-in-N}
There exist positive numbers $\varepsilon_0$ and $\delta_0$ such that
the following conditions hold for any $\varepsilon\in(0,\varepsilon_0]$.
\begin{enumerate}
\item $s(K^\varepsilon)=2n$,
\item $N^\varepsilon\cap\{(x,y,z)\mid x<0\}
\subset N\cap\{(x,y,z)\mid x<0\}$,
\item $b_{\vec v}(K^\varepsilon)=n$, for any unit vector
$\vec v=v_1{\bf i}+v_2{\bf j}+v_3{\bf k}$ with $|v_3|<\delta_0$.
\end{enumerate}
\end{prop}

\noindent{\bf Proof:}
It remains to prove the part (c). As Kuiper did to prove part (a),
we investigate the number of real roots of the function
$t\mapsto(d/dt)K^\varepsilon(t)\cdot\vec v$ for a unit vector
$\vec v=v_1{\bf i}+v_2{\bf j}+v_3{\bf k}$. Approximating $\lambda_1(t)$ and
$\lambda_2(t)$ by finite linear combinations of powers of $\sin t$ and
$\cos t$, we get a curve $\tilde K^\varepsilon$ which is $C^1$-close to
$K^\varepsilon$. We then substitute
$$\cos t=\frac{2w}{1+w^2},\quad \sin t=\frac{1-w^2}{1+w^2}$$
to have
$$\frac d{dt}\tilde K^\varepsilon(t)\cdot\vec v
=\frac{A^{2n}(w)}{(1+w^2)^n}+v_3\cdot\frac{B^{4n}(w)}{(1+w^2)^{2n}}
+\varepsilon\cdot\frac{C^{2N}(w)}{(1+w^2)^{N}}$$
where $A^{2n}$, $B^{4n}$ and $C^{2N}$ are polynomials of degree
$2n$, $4n$ and $2N$, respectively, for some possibly large $N$.
The real roots of this function are the same as those of the polynomial
$$P(w)=A^{2n}(w)\cdot(1+w^2)^{N-n}+v_3\cdot B^{4n}(w)\cdot(1+w^2)^{N-2n}
+\varepsilon\cdot C^{2N}(w).$$
Since
$ A^{2n}(w) =-n\,v_1\sin nt+n\,v_2\cos nt
    =-n\,(v_1^2+v_2^2)^{1/2}\sin(nt-\alpha),$
it has $2n$ real roots.
If $\varepsilon=v_3=0$, they are the real roots of
$P(w)=A^{2n}(w)\cdot(1+w^2)^{N-n},$
each of which is at least one unit away from the remaining
roots $\pm\sqrt{-1}$ of multiplicity $N-n$.
Since the roots of $P(w)$ depend continuously on
$\varepsilon$ and $v_3$, $P(w)$ has exactly $2n$ real roots,
when $\varepsilon$ and $v_3$ are sufficiently small.
One half of them correspond to the local maxima of the projection
$\tilde K^\varepsilon\to\mathbb R\vec v$ and the other half to local minima.
Since $K^\varepsilon$ is $C^1$-close to $\tilde K^\varepsilon$,
part (c) is proved.\qed

\section{Deformations of knots}
In this section, we describe two kinds of deformations
which do not increase the superbridge number.
One is a local deformation and the other is a global one.
\begin{lem}\label{lem:straighten}
Given a knot $K$, let $\bar K$ be a knot obtained by replacing a subarc
of $K$ with a straight line segment joining the end points of the subarc. Then
$s(K)\ge s(\bar K)$.
\end{lem}

\noindent{\bf Proof:} Given a unit vector $\vec v$, let
$g\colon(-1,2)\to\mathbb R\vec v$ be a parametrization of the
orthogonal projection of an open neighborhood of the subarc into $\mathbb
R\vec v$, where the subarc corresponds to the closed interval $[0,1]$. Then
the projection of a neighborhood of the straight line segment in $\bar K$
can be parametrized by
$$\bar g(t)=\begin{cases}(1-t)g(0)+t g(1)&\mbox{if } t\in[0,1]\\
                          g(t)&\mbox{if }t\in(-1,0]\cup[1,2).
              \end{cases}
$$
Since $\bar g$ has no more local maxima than $g$,  we have
$b_{\vec v}(K)\ge b_{\vec v}(\bar K)$ for any $\vec v$.
Therefore $s(K)\ge s(\bar K)$.\qed

\medskip
For a unit vector $\vec v$ and a non-singular linear transformation\break
$\phi\colon\mathbb R^3\to\mathbb R^3$, let $\vec
v^\phi$ denote the unit vector contained in the one-dimensional subspace
$(\phi(\vec v^\perp))^\perp$ satisfying
$\phi(\vec v)\cdot \vec v^\phi>0$. For any subset $A\subset S^2$, we define
$$A^\phi=\{\vec v^\phi\mid \vec v\in A\}.$$

\begin{lem}\label{lem:b_vphi}
Given a unit vector $\vec v\in\mathbb R^3$ and a nonsingular linear
transformation $\phi$ of\/ $\mathbb R^3$, the formulas
\begin{eqnarray*}
&b_{\vec v^\phi}(\phi(K))=b_{\vec v}(K)\\
&b_{\vec v^\phi}(\phi(K)\mid\phi(S))=b_{\vec v}(K\mid S)
\end{eqnarray*}
hold for any knot $K$ and any open subarc $S$ of $K$.
\end{lem}

\noindent{\bf Proof:} At each local maximum point $P$ of the projection
$S\to\mathbb R\vec v$, there is an open disk $d_P$ perpendicular to $\vec v$
and tangent to $S$ at $P$. Then $\phi(d_P)$ is tangent to $\phi(S)$ at $\phi(P)$
and is perpendicular to $\vec v^\phi$.
By the definition of $\vec v^\phi$, $\phi(P)$ is a local maximum point of
the projection $\phi(S)\to\mathbb R\vec v^\phi$ and hence
$b_{\vec v}(K\mid S)\le b_{\vec v^\phi}(\phi(K)\mid \phi(S))$.
Since $(\vec v^\phi)^{\phi^{-1}}=\vec v^{(\phi^{-1}\phi)}=\vec v$,
we also get
$$b_{\vec v}(K\mid S)=
b_{(\vec v^\phi)^{\phi^{-1}}}(\phi^{-1}(\phi(K))\mid\phi^{-1}(\phi(S)))
\ge b_{\vec v^\phi}(\phi(K)\mid \phi(S)).$$
This proves the second formula. Setting $S=K$,
the first formula is obtained.\qed

\medskip
The next proposition easily follows from Lemma~\ref{lem:b_vphi}.
\begin{prop}\label{prop:b_vphi}
Given a knot $K$ and a nonsingular linear transformation $\phi$ of\/ $\mathbb
R^3$, we have $s(\phi(K))=s(K)$. In particular,
if a knot $K$ and a unit vector $\vec v$ satisfy
$b_{\vec v}(K)=s(K)=s[K]$, then $b_{\vec
v^\phi}(\phi(K))=s(\phi(K))=s[K]$.
\end{prop}

\section{Proofs}\label{sec:proofs}
For any $\lambda$ with $0<\lambda\le1$, let
$\phi_\lambda$, $\psi_\lambda$, $\psi$ be the autohomeomorphisms of
$\mathbb R^3$ defined by
\begin{eqnarray*}
\phi_\lambda(x,y,z)&=&(x,y,\lambda z)\\
\psi_\lambda(x,y,z)&=&(1+\lambda-\lambda z,-y,1+\lambda-x)\\
\psi(x,y,z)&=&(-z,-y,-x).
\end{eqnarray*}
The map $\psi$ is the $180^\circ$ rotations about the line
$\{(x,0,z)\mid x+z=0\}$ and the map $\psi_\lambda$ is the composite map
$\phi_\lambda$ followed by the $180^\circ$ rotations about the line
$\{(x,0,z)\mid x+z=1+\lambda\}$.

For any locally one-to-one closed parametrized path
$\gamma\colon S^1\to\mathbb R^3$, we extend the definition of the crookedness
$b_{\vec v}(\gamma)$
by considering the parametrized projection
$t\mapsto\gamma(t)\cdot\vec v:S^1\to\mathbb R\vec v$
instead of the projection $\gamma(S^1)\to\mathbb R\vec v$.
In this way we can consider the crookedness for
finite-fold coverings of knots and singular knots.

\subsection*{Proof of Theorem~\ref{thm:braid-sum}}
Throughout this proof, $\lambda$ is a constant satisfying $0<\lambda\le1/4$,
$\vec v=v_1{\bf i}+v_2{\bf j}+v_3{\bf k}$ is a unit vector, and
$i=1$ or $2$.
We may assume that the knot $K_i$ is parametrized by
$$K_i(t)=((1+\lambda_{i1}(t))\cos n_it,
          (1+\lambda_{i1}(t))\sin n_it,
          \lambda_{i2}(t)+\cos^2n_it)$$
where $\lambda_{i1}$ and $\lambda_{i2}$ are smooth periodic functions
with period $2\pi$ satisfying the conditions corresponding to
(\ref{eqn:in-torus}), (\ref{eqn:as-eta}) and (\ref{eqn:like-eta}), for $i=1,2$.
For any $\varepsilon$ with $0\le\varepsilon\le1$, we define
$K^\varepsilon_1$ and $K^\varepsilon_2$ as in (\ref{eqn:Ke}).
Then $K^\varepsilon_i$ is a knot isotopic to $K_i$ and is the closure of the
$n_i$-braid $K^\varepsilon_i\cap\{(x,y,z)\mid x\le y\le-x\}$
when $0<\varepsilon\le1$. When $\varepsilon=0$, $K^\varepsilon_i$ is an
$n_i$-fold covering of $\eta$.
Since the two knots $\phi_\lambda(K^\varepsilon_1)$ and
$\psi_\lambda(K^\varepsilon_2)$ are tangent at the point
$(1,0,\lambda)$, their union $K_\lambda$ can be regarded as a singular
knot parametrized by
\begin{equation}\nonumber
K_\lambda(t)=
\begin{cases}
\phi_\lambda(K^\varepsilon_1(2t))&\text{if\ }-\pi\le t\le0\\
\psi_\lambda(K^\varepsilon_2(-2t))&\text{if\ }0\le t\le\pi.
\end{cases}
\end{equation}
Then
$\bar K_\lambda=
(K_\lambda-\phi_\lambda(\eta_+)\cup\psi_\lambda(\eta_+))
    \cup S_+\cup S_-$
is a singular knot with only one singular point at
$((1+\lambda)/2,0,(1+\lambda)/2)$ where
$$S_\pm=\{(0,\mp1,0)+s(1+\lambda,\pm2,1+\lambda)\mid 0<s<1\}.$$
By Lemma~\ref{lem:straighten},
$b_{\vec v}(\bar K_\lambda)\le b_{\vec v}(K_\lambda)$.
Since $(1,0,\lambda)$ is a local maximum point of the parametrized projection
$t\mapsto K_\lambda(t)\cdot\vec v$ only if both of the projections
$\phi_\lambda(K^\varepsilon_1)\to\mathbb R\vec v$
and $\psi_\lambda(K^\varepsilon_2)\to\mathbb R\vec v$ have
local maximum at $(1,0,\lambda)$, we have
\begin{equation}\label{eqn:K1+K2}
b_{\vec v}(\bar K_\lambda)\le b_{\vec v}(\phi_\lambda(K^\varepsilon_1))+b_{\vec v}(\psi_\lambda(K^\varepsilon_2)).
\end{equation}
The vectors $\vec w_\pm=\pm(1+\lambda)({\bf i}+{\bf k})+2{\bf j}$,
are parallel to the segments $S_\pm$, respectively.
A computation shows that
\begin{eqnarray*}
\vec w_+\cdot\vec v=
\phantom{-}(1+\lambda)(v_1+v_3)+2v_2\ge\phantom{-}(10-\sqrt{89})/\sqrt{80}>0,
\footnotemark\\
\vec w_-\cdot\vec v= -(1+\lambda)(v_1+v_3)+2v_2\le-(10-\sqrt{89})/\sqrt{80}<0,
\footnotemark
\end{eqnarray*}
whenever $v_3\ge(4\lambda^2+1)^{-1/2}$.%
\addtocounter{footnote}{-1}\footnotetext{%
The equality holds when $\lambda=1/4$ and
$\vec v=-\sqrt{5/89}{\bf i}-8\sqrt{445}{\bf j}+2\sqrt5{\bf k}$.}%
\addtocounter{footnote}{1}\footnotetext{%
The equality holds when $\lambda=1/4$ and
$\vec v=-\sqrt{5/89}{\bf i}+8\sqrt{445}{\bf j}+2\sqrt5{\bf k}$.}
Therefore there exists a number $\delta\in(1/\sqrt2,(4\lambda^2+1)^{-1/2})$
such that
$\vec w_-\cdot\vec v<0<\vec w_+\cdot\vec v\text{\ whenever\ }
v_3\ge\delta$.
At the endpoints $(1+\lambda,\pm1,1+\lambda)$ of $S_\pm$, we have
\begin{eqnarray*}
\lim_{t\to\frac{4n_2-1}{4n_2}\pi^{-}}\frac d{dt}K_\lambda(t)\cdot\vec v
=-2n_2v_3<0<2n_2v_3=
\lim_{t\to\frac\pi{4n_2}^{+}}\frac d{dt}K_\lambda(t)\cdot\vec v
\end{eqnarray*}
if $v_3>0$.
Therefore there exist open arcs $\tilde S_\pm$ of $\bar K_\lambda$,
containing the closures of $S_\pm$, respectively,
satisfying
$b_{\vec v}(\bar K_\lambda\mid \tilde S_+\cup\tilde S_-)=0$
whenever $v_3\ge\delta$.
Similarly we also have
$b_{\vec v}(\bar K_\lambda\mid\tilde S_+\cup\tilde S_-)=0$
whenever $v_1\le-\delta$.

By Lemma~\ref{lem:1-or-0}, Lemma~\ref{lem:b_vphi}, and
the last two conditions, we have\footnote{$\overline\eta_+$ is the closure of
$\eta_+$.}
\begin{eqnarray}\label{eqn:K1+K2-1}
\nonumber b_{\vec v}(\bar K_\lambda)
&\le&b_{\vec v}(\phi_\lambda(K^\varepsilon_1)\mid
            \phi_\lambda(K^\varepsilon_1-\overline\eta_+))
    +b_{\vec v}(\bar K_\lambda\mid \tilde S_+\cup\tilde S_-)\\
&&\nonumber\phantom{b_{\vec v}(\phi_\lambda(K^\varepsilon_1))}
    +b_{\vec v}(\psi_\lambda(K^\varepsilon_2)\mid
            \psi_\lambda(K^\varepsilon_2-\overline\eta_+))\\
&\le& b_{\vec v}(\phi_\lambda(K^\varepsilon_1))+b_{\vec v}(\psi_\lambda(K^\varepsilon_2))-1
\end{eqnarray}
whenever
$\vec v\in N^{\phi_\lambda}\cup Q_\delta\cup\psi(N^{\phi_\lambda}\cup Q_\delta)$
where $Q_\delta=\{(x,y,z)\in S^2\mid x>0,z>\delta\}$.
By Proposition~\ref{prop:Ne-in-N}~(a)--(b),
we may assume that
\begin{eqnarray}
\label{eqn:sKi=2ni}&s(K^\varepsilon_i)=2n_i\\
\label{eqn:Ne-in-N}
&(N^\varepsilon_i)^{\phi_\lambda}\subset N^{\phi_\lambda}\cup Q_\delta
\end{eqnarray}
where $N^\varepsilon_i=\{\vec v\in S^2\mid
      v_3>0,b_{\vec v}(K^\varepsilon_i)=2n_i\}.$
Since
\begin{eqnarray*}
&&\vec v^{\phi_\lambda}\cdot{\bf k}=v_3(\lambda^2(1-v_3^2)+v_3^2)^{-1/2},\\
&&\vec v^{\psi_\lambda}\cdot{\bf i}=-v_1(\lambda^2(1-v_1^2)+v_1^2)^{-1/2},
\end{eqnarray*}
Proposition~\ref{prop:Ne-in-N}~(c) implies that
\begin{equation}
\label{eqn:Kei-ni}
\begin{array}{l}
\displaystyle b_{\vec v}(\phi_\lambda(K^\varepsilon_1))=n_1
\text{\ whenever\ } |v_3|\le1/\sqrt2\\
\displaystyle b_{\vec v}(\psi_\lambda(K^\varepsilon_2))=n_2
\text{\ whenever\ } |v_1|\le1/\sqrt2
\end{array}
\end{equation}
provided $\lambda$ is sufficiently small.
By (\ref{eqn:K1+K2}) and (\ref{eqn:Kei-ni}), we get
\begin{equation}\label{eqn:bK'-1}
b_{\vec v}(\bar K_\lambda)\le n_1+n_2
\text{\ if\ } \max\{|v_1|,|v_3|\}\le1/\sqrt2.
\end{equation}
By (\ref{eqn:K1+K2}), (\ref{eqn:sKi=2ni}) and (\ref{eqn:Kei-ni}), we get
\begin{equation}\label{eqn:bK'-2}
b_{\vec v}(\bar K_\lambda)\le
\begin{cases}
2n_1+n_2-1&
\text{if\ } \pm\vec v\notin(N^\varepsilon_1)^{\phi_\lambda},
            |v_3|>1/\sqrt2\\
n_1+2n_2-1&
\text{if\ }\pm\vec v\notin\psi((N^\varepsilon_2)^{\phi_\lambda}),
            |v_1|>1/\sqrt2
\end{cases}
\end{equation}
By (\ref{eqn:K1+K2-1}), (\ref{eqn:sKi=2ni}), (\ref{eqn:Ne-in-N})
and (\ref{eqn:Kei-ni}), we get
\begin{equation}\label{eqn:bK'-3}
b_{\vec v}(\bar K_\lambda)\le
\begin{cases}
2n_1+n_2-1&
\text{if\ } \pm\vec v\in(N^\varepsilon_1)^{\phi_\lambda}\cup Q_\delta\\
n_1+2n_2-1&
\text{if\ }\pm\vec v\in\psi((N^\varepsilon_2)^{\phi_\lambda}\cup Q_\delta)
\end{cases}
\end{equation}
For the last two formulas, we used the fact
$b_{-\vec v}(\bar K_\lambda)=b_{\vec v}(\bar K_\lambda)$.
For a very small positive number $\epsilon$, let
$\bar S_+=S_+\cup\{(\cos t,\sin t,\lambda\cos^2t)\mid -\pi/2-\epsilon\le t\le-\pi/2\}$ and let $\check S_+$ be the line segment joining the endpoints of
$\bar S_+$. By the conditions (\ref{eqn:in-torus}),
(\ref{eqn:as-eta}) and (\ref{eqn:like-eta}), the knot
$\check K_\lambda=(\bar K_\lambda-\bar S_+)\cup\check S_+$
is a knot representing $K_1\sharp K_2$. By Lemma~\ref{lem:straighten},
(\ref{eqn:bK'-1}), (\ref{eqn:bK'-2}) and (\ref{eqn:bK'-3}),
we have
\begin{equation}\nonumber
b_{\vec v}(\check K_\lambda)\le b_{\vec v}(\bar K_\lambda)
\le\max\{2n_1+n_2,n_1+2n_2\}-1.\kern1.2cm\qed\kern-1.2cm
\end{equation}

\subsection*{Proof of Theorem~\ref{thm:torus-sum}}
Let $K_i$ be a torus knot of type $(p_i,q_i)$ where $p_i$ and $q_i$ are
coprime integers satisfying $2\le p_i<q_i$, for $i=1,2$.
This proof breaks into three cases.

\noindent{\sc Case 1.}
Suppose that the inequality $2\le p_i<q_i/2$ holds for $i=1,2$.
In this case, we have $\beta[K_i]=b[K_i]=s[K_i]/2=p_i$.
Therefore a direct application of Theorem~\ref{thm:braid-sum} shows that
$s[K_1\sharp K_2]\le \max\{2p_1+p_2,p_1+2p_2\}-1$.

\noindent{\sc Case 2.}
Suppose that the inequality $2\le p_i<q_i<2p_i$ holds for $i=1,2$.
In this case, $\beta[K_i]=b[K_i]=p_i$ and $s[K_i]=q_i$.
As shown in~\cite{jin-poly},
$K_i$ can be represented by a polygonal knot
$\tau_i=\tau_i(\alpha_i)$
of $2q_i$ edges embedded on the torus $H_{\alpha_i}\cup H_{\beta_i}$
where
$$H_\theta=\{(x,y,z)\mid
x^2+y^2-z^2\sin^2\frac\theta2=\cos^2\frac\theta2,
|z|\le 1\},$$
$\pi p_i/q_i<\alpha_i<2\pi p_i/q_i$ and $\alpha_i+\beta_i=\pi$.
The knot $K_i$ has $2q_i$ vertices; $q_i$ on each of the two unit circles
$\{(x,y,\pm1)\mid x^2+y^2=1\}$.
By (\ref{eq:2s<p}), we know that $s(K_i)=q_i$.
We may assume that $K_i$ has a vertex at $(1,0,1)$.
We define
\begin{eqnarray*}
N_i&=&\{\vec v\in S^2\mid \vec v\cdot{\bf k}>0,b_{\vec v}(K_i)=q_i\},\\
M_i&=&\{\vec v\in S^2\mid \text{The projection $K_i\to\mathbb R\vec v$ has
a local minimum at $(1,0,1)$}\}.
\end{eqnarray*}
For any $\vec v\in N_i$, the $q_i$ vertices of $K_i$ on the circle
$\{(x,y,1)\mid x^2+y^2=1\}$ are local maximum points of the projection
$K_i\to\mathbb R\vec v$.
Let $t\mapsto K_i(t)$ parametrize $K_i$ modulo $2\pi$ with $K_i(0)=(1,0,1)$,
as a closed $p_i$-braid around the $z$-axis. The
singular knot $K_\lambda$ given by the parametrization
$$ K_\lambda(t)=\begin{cases}
\phi_\lambda(K_1(2t))&\text{if\ }-\pi\le t\le0\\
\psi_\lambda(K_2(-2t))&\text{if\ }0\le t\le\pi
\end{cases}$$
has only one singular point at $(1,0,\lambda)$.
Straightening an arc near the singular point, we get a knot representing
$K_1\sharp K_2$ whose crookedness is not bigger than that of $K_\lambda$
in any direction.
As $\lambda$ approaches zero, $N_1$ shrinks to the north pole $(0,0,1)$
whereas $M_1$ approaches a region of positive area containing
the point $(-1,0,0)$. Therefore, for a sufficiently small $\lambda$, we have
\begin{equation}\label{eqn:N-in-M}
(N_1)^{\phi_\lambda}\subset \psi((M_2)^{\phi_\lambda}),\text{\ and\ }
\psi((N_2)^{\phi_\lambda})\subset (M_1)^{\phi_\lambda},
\end{equation}
and as in (\ref{eqn:Kei-ni}), we also have
\begin{equation}
\label{eqn:Ki-ni}
\begin{array}{l}
\displaystyle b_{\vec v}(\phi_\lambda(K_1))=p_1
\text{\ whenever\ } |\vec v\cdot{\bf k}|\le1/\sqrt2,\\
\displaystyle b_{\vec v}(\psi_\lambda(K_2))=p_2
\text{\ whenever\ } |\vec v\cdot{\bf i}|\le1/\sqrt2.
\end{array}
\end{equation}
By (\ref{eqn:N-in-M}), if
$\pm\vec v\in (N_1)^{\phi_\lambda}\cup\psi((N_2)^{\phi_\lambda})$,
the point $(1,0,\lambda)$ in not a local maximum point of $K_\lambda$.
Therefore we have
\begin{equation}\nonumber
b_{\vec v}(K_\lambda)=\begin{cases}
q_1+p_2-1&\text{if\ }\pm\vec v\in (N_1)^{\phi_\lambda}\\
p_1+q_2-1&\text{if\ }\pm\vec v\in \psi((N_2)^{\phi_\lambda}).
\end{cases}
\end{equation}
By (\ref{eqn:N-in-M}) and (\ref{eqn:Ki-ni}), we obtain
$$b_{\vec v}(K_\lambda)\le b_{\vec v}(\phi_\lambda(K_1))
+b_{\vec v}(\psi_\lambda(K_2))<\max\{q_1+p_2,p_1+q_2\},$$
if $\pm\vec v\notin (N_1)^{\phi_\lambda}\cup\psi((N_2)^{\phi_\lambda})$.
Therefore $b_{\vec v}(K_\lambda)\le\max\{q_1+p_2,p_1+q_2\}-1$, for any
unit vector $\vec v$.

\noindent{\sc Case 3.}
Suppose that the inequalities $2\le p_1<q_1/2$ and $2\le p_2<q_2<2p_2$ hold.
Let $K_1\,(=K^\varepsilon_1)$ and $K_2$ be embedded and parametrized
as in the Proof of
Theorem~\ref{thm:braid-sum} and in {\sc Case~2}, respectively.
We consider the singular knot parametrized by
$$ K_\lambda(t)=\begin{cases}
\phi_\lambda(K^\varepsilon_1(2t))&\text{if\ }-\pi\le t\le0\\
\psi_\lambda(K_2(-2t))&\text{if\ }0\le t\le\pi.
\end{cases}$$
We replace the arc $\phi_\lambda(\eta_+)$ of $\phi_\lambda(K^\varepsilon_1)$
by the broken line joining the three points $(0,-1,0)$, $(1,0,\lambda)$
and $(0,1,0)$, consecutively, to get a new singular knot $\bar K_\lambda$.
The remaining argument will be very similar to that of {\sc Case~2}.\qed

\bigskip
{\small
\noindent{\sc Acknowledgement:} This work was done while the author
was visiting the University of British Columbia during the academic year
1998--99. He is grateful to the members of the Department of Mathematics
of UBC, especially Dale Rolfsen with whom he had many helpful discussions.

%

\begin{thebibliography}{XX}
%
\bibitem{burde-zieschang}
    G.\ Burde and H.\ Zieschang,
    \textit{Knots},
    de Gruyter Studies in Mathematics vol.~5,
    Walter de Gruyter, Berlin, New York, 1985
\bibitem{calvo}
    J.\ Calvo,
    Geometric knot theory: the classification of spatial polygons
         with a small number of edges,
    Ph.D.\ thesis, University of California Santa Barbara, 1998
\bibitem{calvo-millett}
    J.\ Calvo and K.C.\ Millett,
    Minimal edge piecewise linear knots,
    \textit{Ideal Knots}
    (Series on Knots and Everything vol.~19, World Scientific,
    1998) 107--128.
\bibitem{doll}
    H.\ Doll,
    A generalized bridge number for links in 3-manifolds,
    Math.\ Ann.\ \textbf{294}(1992) 701--717.
\bibitem{jeon-jin-3}
    C.B.\ Jeon and G.T.\ Jin,
    There are only finitely many 3-superbridge knots,
    J.\ Knot Theory Ramifications
    (Special issue of {\em Knots in Hellas 1998}\/) to appear.
\bibitem{jeon-jin}
    C.B.\ Jeon and G.T.\ Jin,
    A computation of superbridge index of knots,
    in preparation.
\bibitem{jin-poly}
    G.T.\ Jin,
    Polygon indices and superbridge indices of torus knots and links,
    J.\ Knot Theory Ramifications \textbf{6}(1997) 281--289.
\bibitem{jin-kim}
    G.T.\ Jin and H.S.\ Kim,
    Polygonal knots, J.\ Korean Math.\ Soc.\ \textbf{30}(1993) 371--383.
\bibitem{kuiper}
    N.\ Kuiper,
    A new knot invariant,
    Math.\ Ann.\ \textbf{278}(1987) 193--209.
\bibitem{meissen}
    M.\ Meissen,
    Edge number results for piecewise-linear knots,
    \textit{Knot theory}
    (Banach Center Publications vol.~42, Warszawa, 1998)
    235--242.
\bibitem{millett}
    K.C.\ Millett,
    Coordinates of nine stick $9_{44}$ and $9_{46}$,
    Email correspondence.
\bibitem{milnor1}
    J.W.\ Milnor,
    On the total curvature of knots,
    Ann.\ Math.\ \textbf{52}(1950) 248--257.
\bibitem{negami}
    S.\ Negami,
    Ramsey theorems for knots, links and spatial graphs,
    Trans.\ Am.\ Math.\ Soc.\ \textbf{324}(1991), no.\ 2, 527--541.
\bibitem{rolfsen}
    D.\ Rolfsen,
    \textit{Knots and Links}
    Mathematics Lecture Series 7, Publish or Perish, 1976
\bibitem{scharein}
    R.\ Scharein, Minimum Stick Candidates, \newline
    {\tt http://www.cs.ubc.ca/nest/imager/contributions/scharein/sa/msc.html}
\bibitem{schubert}
    H.\ Schubert,
    \"Uber eine numerische Knoteninvariente,
    Math.\ Z.\ \textbf{61}(1954) 245--288.
%
\end{thebibliography}
%

\end{document}